\documentclass[11pt]{article}
\usepackage[textures]{epsfig,psfrag}
\usepackage{amsmath,amssymb,amsthm,amscd,tabularx,array,fancyhdr}
\newfont{\extra}{msbm10 scaled\magstep1}

\newcommand{\extr}[1]{\mbox{\extra #1}}

\newcommand{\sect}[1]{\setcounter{equation}{0}\section{#1}}

\newtheorem{theorem}{Theorem}[section] 
 
\newtheorem{lemma}{Lemma}[section]

\def\bee{\begin{equation}}
\def\eee{\end{equation}}
\def\bea{\begin{eqnarray}}
\def\eea{\end{eqnarray}}

\newcommand{\rad}[1]{ {\text{\rm ad}^{\text{\rm r}}_{#1}}}

\def\C{\extr C}
\def\K{\extr K}
\def\N{\extr N}
\def\R{\extr R}

\newcommand{\lact}{\triangleright}
\newcommand{\ract}{\triangleleft}
\newcommand{\lcact}{\blacktriangleleft}

\begin{document}

\begin{center}{ \LARGE \bf 
Local representations of \\[0.5cm]
the $(1+1)$ quantum extended Galilei algebra}
\end{center}
\vskip0.25cm

\begin{center}
Oscar Arratia $^1$ and Mariano A. del Olmo $^2$
\vskip0.25cm

{ \it $^{1}$ Departamento de  Matem\'atica Aplicada a la  Ingenier\'{\i}a,  \\
\vskip0.10cm
$^{2}$ Departamento de  F\'{\i}sica Te\'orica,\\
\vskip0.10cm
 Universidad de  Valladolid, 
 E-47011, Valladolid,  Spain}

\vskip0.15cm

E. mail: oscarr@wmatem.eis.uva.es, olmo@fta.uva.es
 
\end{center}
 
\vskip1.5cm
\centerline{\today}
\vskip1.5cm

\begin{abstract}
We present the $q$--deformed counterpart of the local representations of the 
$(1+1)$ extended  Galilei group. These representations act on the space of wavefunctions defined in
the space-time.  As in the classical case the $q$--local representations are reducible and the
condition of irreducibility is given by the $q$--Casimir equation that in the limit of the deformation
parameter going to zero becomes the Schr\"odinger equation of a free particle.
\end{abstract}
\newpage

\sect{Introduction}

It is well known that in Quantum Mechanics the symmetry group $G$ of a quantum system
is  realized by local unitary operators $U$ acting on the space of wavefunctions $\psi$ 
defined on the space-time manifold (a homogeneous space of $G$) $X$ in the form \cite{cos85} 
$$
 \psi'(x')\equiv (U(g) \psi)(g\,x) 
= A(g,x) \psi(x), \qquad g\in G,\ x\in X,
$$
being $A(g,x)$ a matrix-valued function. When the  wavefunctions are single component 
$A(g,x)$ is a  phase function $A(g,x)=\exp^{i\zeta (g,x)}$.  
In general, the operators $U(g)$ do not close a true 
representation  but a projective or `up to a factor' representation 
of $G$ called {\sl local realization} 
$$
U(g_2)U(g_1) =\omega(g_2,g_1) U(g_2g_1),\quad g_2,g_1 \in G.
\label{UU}
$$
The function $\omega :G\times G \to U(1)$ is the factor system of the 
realization and it is a 2--cocycle. 
To construct  the local realizations of a given group $G$ the  first step is to linearize the problem. 
In other words,  instead of computing directly  the representations up to a factor of $G$ we can
obtain them from the linear local  representations of a new group $\bar G$, which is
a central extension of  $G$ by an Abelian group
$A$. It can be shown that $A$ is related with the second
cohomology group of $G$ \cite{cos84}. So, the computation of ${\bf H}^2(G,U(1))$ can be done by
solving the equivalent problem of the central extensions of $G$ by $U(1)$. 
It can be proved \cite{cos85} that the local representations can be obtained by induction from the
finite dimensional representations of the subgroup $\bar\Gamma$ of $\bar G$ such that  
${\bar G}/{\bar\Gamma}\simeq X$. 

On the other hand, the study of the $q$--deformed symmetries of the space-time is one of the 
main applications of  quantum kinematical groups in order to establish the structure of the space-time 
at the level of the distances and energies of Planck scale. This is an important problem to be solved in
relation with the  statement of a quantum gravity theory.

So, it looks interesting to study the $q$--deformed version of the local representations of kinematical
groups in order to construct  models of deformed quantum mechanics, where the deformation parameter
of the quantum algebra is introduced as a fundamental length that can be interpreted as a lattice
spacing. Attempts in this direction has been presented in Ref.~\cite{italiani92}, \cite{Bon92} and
\cite{luki}. Also recent progress has been made in $q$-deformed quantum field theory  by Wess
and coworkers \cite{wess}. This work can be seen from the point of view of non-commutative geometry 
\cite{connes,varilly01}.

In \cite{olmo00} we constructed the
induced representations of some quantum versions of the Galilei algebra in $(1+1)$--dimensions, among
them  the quantum extended Galilei algebra
$U_q(\overline{\mathfrak{g}(1,1)})$. However, the  construction of induced representations of ``local
type"  for $U_q(\overline{\mathfrak{g}(1,1)})$ presented two problems: the generators associated to the
space-time coordinates do not appeared isolated on a side of the monomial basis chosen, and the linear
space spanned by the elements related with $\bar\Gamma$ does not close a subalgebra of
$U_q(\overline{\mathfrak{g}(1,1)})$. 
In this paper we have solved these difficulties and we present the `local' induced representations of
$U_q(\overline{\mathfrak{g}(1,1)})$. The solution proposed by us does without the requirement of duality
 or pairing between algebras replacing it with a less restrictive condition of duality between linear
spaces.  The price payed is to drop the coalgebra sector.

The organization of the paper is as follows. Section~\ref{q-representations} is devoted to present a
brief review about the theory of induced representations of quantum algebras developed in
Ref.~\cite{olmo00}--\cite{olmo01}. In Section \ref{q-galilei} we describe the quantum extended  Galilei
algebra in $(1+1)$--dimensions.  The main results of the paper about the construction of the local
representations of the quantum  algebra mentioned above are presented in Section
\ref{q-represgalilei}. We finish with some conclusions and remarks.

\sect{Induced representations of quantum algebras}\label{q-representations}

The main results about the theory of induced representation of quantum  algebras
will be summarized in this section (see \cite{olmo00} and \cite{olmo01}   for a
detailed description).

We will use  multi-index notation: let $A$ be an
algebra generated by the  elements $(a_1, a_2, \ldots, a_r)$ such that the ordered monomials 
$a^n := a_1^{n_1} a_2^{n_2} \cdots a_r^{n_r} \in A, \  
n=(n_1, n_2, \ldots, n_r) \in \mathbb{N}^r$,
form a basis of the linear space  underlying to $A$.
An arbitrary product of generators of $A$ is written in a normal ordering if it is 
expressed in terms of the (ordered) basis  $(a^n)_{n \in \mathbb{N}^r}$.  The unit of $A$,
$1_A$, is denoted by  $a^0 $ ($0 \in \N^n$). Multi-factorials
and multi-deltas are defined by
$l! = \prod_{i=1}^n l_i!$ and $\delta_l^m= \prod_{i=1}^n \delta_{l_i}^{m_i}$.

We will denote the different actions involved by the following symbols (or their symmetric
for the corresponding right actions and  coactions):
$\lact$  ($\lcact$) (actions  (coactions)),
$\vdash$ (induced and inducting representations) and
$\succ$   (regular actions).

Let $H$ be a Hopf algebra and $V$ a linear vector space over a field
$\K$ ($\R$ or $\C$). The triplet $(V, \triangleright, H)$ 
is said to be a left $H$--module if $\alpha$ is a left action of $H$ on $V$, i.e., 
 a linear map $\alpha: H \otimes V \to V$
($ \alpha : (h \otimes v)  \mapsto\alpha(h \otimes v)\equiv h \triangleright v$) 
such that 
$$
 h_1\triangleright (h_2\triangleright v)= (h_1h_2) \triangleright v, \quad
 1_H \triangleright v = v,
\qquad \forall h_1,h_2 \in H,\ \forall v \in V.
$$
Right $H$--modules can be defined in a similar way.

There are two canonical modules associated to any pair of Hopf algebras, $H$ and $H'$,  related by a
nondegenerate pairing $ \langle \, \cdot \, , \, \cdot \, \rangle$, in the sense of being
simultaneously  left and right nondegenerate pairing \cite{charipressley} 
(in this case $(H, H', \langle \, \cdot \, , \, \cdot \, \rangle)$ 
will be called a nondegenerate triplet):

\ i) The left regular module $(H, \succ, H )$ with  action  
$$
 h_1 \succ h_2 = h_1h_2,\qquad  \forall  h_1,h_2 \in H.
$$

ii) The right coregular module  $(H', \prec, H )$ with  action  defined by
$$
\langle h_2, h' \prec h_1 \rangle =
\langle h_1 \succ h_2, h' \rangle,
\qquad \forall h_1,h_2 \in  H, \quad \forall h' \in H'.
$$
Using the coproduct ($\Delta(h')= h'_{(1)}\otimes h'_{(2)}$) in $H'$ the action takes the form:
$h' \prec h = \langle h, h'_{(1)} \rangle h'_{(2)}$.

The induction and coinduction algorithms of algebra representations are adapted to the
Hopf algebras as follows.
 Let $(H, H', \langle \, \cdot \, , \, \cdot \, \rangle)$ be a nondegenerate triplet 
 and  $(V, \triangleright ,K)$  a left $K$--module with $K$ a subalgebra of $H$. The
carrier space, $\K^\uparrow$, of the coinduced representation is the subspace of
$H'\otimes V$ with elements $f$ verifying
 \begin{equation}\label{coinducidar}
 \langle f, k h \rangle = k \triangleright  \langle f, h \rangle, \qquad
\forall k \in K, \ \forall h \in H.
\end{equation}
The pairing  used in expression (\ref{coinducidar}) is $V$--valued and is defined  by
$\langle h' \otimes v, h \rangle= \langle h', h \rangle v$, with
$h \in H,\  h' \in H',\  v \in V$.
The action $h\triangleright f$ on the coinduced module is determined by
$$
  \langle h_1 \triangleright f, h_2 \rangle =
\langle f, h_2h_1 \rangle, \qquad \forall h_2 \in H.
$$

Let $(\K, \triangleright, K)$ be a  coinducing module.
The carrier space of the coinduced representation is the subspace
of $H' \otimes \K \simeq H'$ composed by the elements $\varphi$
that verify the equivariance condition
 $\varphi \prec k = (1 \dashv k)\varphi ,\ \forall k \in K.$
The action of $H$ on $\K^\uparrow$ induced by the action of $K$ on $\K$ is given by
$$
\langle h_2 \triangleright \varphi, h_1 \rangle = \langle \varphi, h_1h_2\rangle, 
\qquad \forall h_1, h_2 \in H, \quad \forall \varphi \in \K^\uparrow ,
$$
or explicitly by 
$ h \triangleright \varphi\equiv h \succ \varphi =  \langle h, \varphi_{(2)}\rangle
\varphi_{(1)}$.

Note that to describe the induced module  the
right, $({H'}, \prec, H)$, and left, $(H', \succ, H )$,  coregular modules are both
pertinent, the former  to determine the carrier space and the last to obtain the
induced action.

Let $(H, {H'}, \langle \, \cdot \, , \, \cdot \,
\rangle)$ be a nondegenerate triplet with finite sets of generators,
$\{h_1,\ldots, h_n\}$ of $H$ and $\{\varphi^1, \ldots, \varphi^n\}$ of $H'$, such that 
$\{ h_l=h_1^{l_1}\cdots h_n^{l_n}\}_{l \in \N^n}$ and
$\{ \varphi^m= (\varphi^1)^{m_1} \cdots (\varphi^n)^{m_n}\}_{m \in \N^n}$ ($l=(l_1,
\ldots, l_n)$, $m=(m_1, \ldots, m_n)$) are bases of
$H$ and $H'$, respectively. 
The action on the coregular module
$({H'}, \succ, H)$ is obtained after to compute the action of the  generators
$$
h_i \succ \varphi^j= \sum_{k \in \N^n} \alpha_{ik}^j \varphi^k,
\qquad i, j \in \{  1, 2, \ldots, n \},
$$
and to extend it to the ordered polynomial
$\varphi^j= (\varphi^1)^{j_1}\cdots (\varphi^n)^{j_n}$ using the
compatibility relation between the action
 and the algebra structure in $H'$. Thus, we get
\begin{equation}\label{primera}
 h \succ (\varphi \psi)= (h_{(1)} \succ \varphi)
(h_{(2)} \succ \psi),\qquad
 h \succ 1_{H'} = \epsilon(h) 1_{H'}.
 \end{equation}
In order  to write  explicitly the expression of the action on
a generic ordered polynomial we take into account that:

1) There is a natural
representation $\rho$ of  $H$  associated to  $(H,\prec,H)$   
$$
 [\rho(h_2)](h_1)= h_1 \prec h_2.
$$

2) The action on $({H'},\succ,H)$ can be expressed in terms of 
$\rho$  using the adjoint with respect to the pairing
($f^\dagger :H'\to H'$ is
the adjoint of $f :H\to H$ if $\langle \, h ,f^\dagger(h') \rangle= \langle \, f(h)
,h'\rangle$) defined by
\begin{equation}\label{leftcoregularaction}
 h \succ \varphi= [\rho(h)]^\dagger(\varphi).
\end{equation}

If the bases $\{ h_l\}_{l \in \N^n} $ and $\{ \varphi^m\}_{m \in \N^n}$  are dual (i.e.
$\langle h_l, \varphi^m \rangle = l! \; \delta_l^m, \  \forall l, m \in \N^n$) 
we  define  `multiplication' operators $\overline{h}_i$,
$\overline{\varphi}^j$ and formal derivatives
${\partial}/{\partial h_i}$,
${\partial}/{\partial \varphi^j}$ by
$$
 \begin{array}{l}
  \overline{h}_i (h_1^{l_1} \cdots h_i^{l_i} \cdots h_n^{l_n})=
  h_1^{l_1} \cdots h_i^{l_i+1} \cdots h_n^{l_n}, \\[2mm]
 \overline{\varphi}_i \left( (\varphi^1)^{m_1} \cdots (\varphi^i)^{m_i} \cdots 
(\varphi^n)^{m_n}\right) = (\varphi^1)^{m_1} \cdots (\varphi^i)^{m_i+1} \cdots 
(\varphi^n)^{m_n},\\[2mm]
\displaystyle{\frac{\partial}{\partial h_i}} (h_1^{l_1} \cdots h_i^{l_i} \cdots h_n^{l_n}) =
 l_i \; h_1^{l_1} \cdots h_i^{l_i-1} \cdots h_n^{l_n}, \\[2mm]
\displaystyle {\frac{\partial}{\partial \varphi^i}}
\left( (\varphi^1)^{m_1} \cdots (\varphi^i)^{m_i} \cdots (\varphi^n)^{m_n}\right) 
= m_i \; (\varphi^1)^{m_1} \cdots (\varphi^i)^{m_i-1} \cdots (\varphi^n)^{m_n}.
\end{array}
$$
The adjoint operators are given by 
$\overline{h}_i^\dagger=  {\partial}/{\partial \varphi^i}$ and 
$\overline{\varphi}^{i\dagger} =  {\partial}/{\partial h_i}$.

\sect{Quantum extended Galilei algebra $U_q(\overline{\mathfrak{g}(1,1)})$}
\label{q-galilei}

The quantum extended Galilei group, $F_q(\overline{G(1,1)})$, can be obtained as a
deformation of a  Lie--Poisson structure defined by a nontrivial  1--cocycle of the extended
Galilei group \cite{olmo91}, $\overline{G(1,1)}$, with values on $\Lambda^2 \overline{\mathfrak{g}(1,1)}$
\cite{bgst98}.  Thus, $F_q(\overline{G(1,1)})$ is the Hopf algebra generated by $\mu, x, t$
and $v$ with nonvanishing commutation relations 
$$
[\mu,x]=-2a\mu, \qquad [\mu,v]= a v^2, \qquad [x,v]= 2av.
$$
The coproduct, counit and antipode are
$$
\begin{array}{c}
\begin{array}{ll}
\Delta \mu= \mu \otimes 1 + 1 \otimes \mu + v \otimes x + 
\frac{1}{2} v^2 \otimes t, &\qquad
\Delta x= x \otimes 1 + 1 \otimes x + v \otimes t, \\[2mm]
\Delta t = t \otimes 1 + 1 \otimes t ,&\qquad
\Delta v = v \otimes 1 + 1 \otimes v;
 \end{array} \\ \\[-2.5mm]
\epsilon(\mu)=\epsilon(x)=\epsilon(t)=\epsilon(v)=0 ; \\[2mm]
S(\mu)= -\mu + x - \frac{1}{2} v^2 t, \quad S(x)= -x +t v , \quad
S(t)= -t, \quad S(v)= -v.
\end{array}
$$

The dual (enveloping) algebra $U_q(\overline{\mathfrak{g}{(1,1)}})$
\cite{Bon92} is generated by $I, P, H$ and $N$, in such a way that its
pairing with $F_q(\overline{G(1,1)})$ is given by
\begin{equation}\label{pairing}
\langle I^p P^q H^r N^s, \mu^{p'} x^{q'} t^{r'} v^{s'}\rangle =
p! q! r! s! \; \delta^p_{p'} \delta^q_{q'} \delta^r_{r'} \delta^s_{s'}.
\end{equation}
The duality relations fix the Hopf algebra structure in
$U_q(\overline{\mathfrak{g}{(1,1)}})$. The nonvanishing commutation relations are
\begin{equation}\label{conmutadores}
  [I,N]=-a e^{-2aP} I^2, \qquad  [P,N]= - e^{-2a P} I, 
\qquad  [H, N]=- {\displaystyle \frac{1- e^{-2 a P}}{2a}}.
 \end{equation}
Note that  $I,P$ and  $H$ close a commutative subalgebra
but  now $I$ is not central. However, the commutation relations (\ref{conmutadores}) are 
simpler in terms of the  generators $(M,P,H,K)$, whose nonzero Lie brackets are 
\begin{equation} \label{relconmugext}
 [P,K]=-M, \qquad [H,K]=- {\displaystyle \frac{\sinh(aP)}{a}}.
\end{equation}
 where $M= e^{-a P} I $ and  $K= e^{aP} N$. However, the pairing has a 
`diagonal' expression using $(I,P,H,N)$. 
 
The  coproduct, counit and antipode can be written as
$$
\begin{array}{c}
\begin{array}{ll}
\Delta M= M \otimes e^{-a P} + e^{a P} \otimes M, &\qquad
\Delta P= P \otimes 1 + 1 \otimes P, \\[2mm]
\Delta H= H \otimes 1 + 1 \otimes H, &\qquad
\Delta K= K \otimes e^{-a P} + e^{a P} \otimes K;
\end{array} \\ \\[-2.5mm]
\epsilon(M)= \epsilon(P)=\epsilon(H)=\epsilon(K)=0; \\[2mm]
S(M)= -M, \quad S(P)= -P, \quad S(H)= -H,
\quad S(K)= -K- aM;
\end{array}
$$

From the commutation relations (\ref{relconmugext})
we see that $M,P$ and $H$ span a commutative subalgebra where the generator $K$
acts like a derivation by means of the  commutator. Moreover, since  $M$ is a central
generator we will identify it with the mass operator of a physical system.

\sect{Local representations of $U_q(\overline{\mathfrak{g}(1,1)})$}
\label{q-represgalilei}

We mentioned before in the Introduction the difficulties that prevented the  construction of local
representations  for $U_q(\overline{\mathfrak{g}(1,1)})$ in \cite{olmo00}, related with that the
generators $x$ and $t$ do not appear isolated on a side of the monomial basis, and with that the linear
space spanned by the elements $I^m N^n$ does not close a subalgebra of
$U_q(\overline{\mathfrak{g}(1,1)})$.  Here we will use an alternative basis of
$U_q(\overline{\mathfrak{g}(1,1)})$ determined by $(M,P,H,K)$. However, new troubles appear because 
the pairing between $U_q(\overline{\mathfrak{g}(1,1)})$ and
$F_q(\overline{G(1,1)})$ is not like (\ref{pairing})  (i.e., 
$\langle M^p P^q H^r K^s, \mu^{p'} x^{q'} t^{r'} v^{s'}\rangle \neq
p! q! r! s! \; \delta^p_{p'} \delta^q_{q'} \delta^r_{r'} \delta^s_{s'}$). 
But in order  to avoid them 
let us consider a unital associative algebra 
$\cal A$ spanned by four elements
$v,\mu,x, t$, whose commutation relations are not specified. Although, let us suppose 
that the ordered monomials 
$v^m \mu^n x^p t^q$ determine a basis of the underlying linear space of $\cal A$, and
the nondegenerate bilinear form defined on 
$U_q(\overline{\mathfrak{g}(1,1)})\times {\cal A}$ is
\begin{equation} \label{pairinggext2}
\langle K^m M^n P^p H^q, v^{m'} \mu^{n'} x^{p'} t^{q'} \rangle = m! n! p! q! \; 
\delta^m_{m'} \delta^n_{n'} \delta^p_{p'} \delta^q_{q'}.
 \end{equation}
So, in the following we will work only with the algebra sector
of $U_q(\overline{\mathfrak{g}(1,1)})$. 
Note that the nondegenerate pairings play a fundamental role to transfer the module 
structure to the dual spaces.

The first step is to know the structure of the module-coalgebra  
$(U_q(\overline{\mathfrak{g}(1,1)}), \prec, U_q(\overline{\mathfrak{g}(1,1)}))$. 
Taking into account that
\begin{equation} 
M^n P^p H^q K= K M^n P^p H^q - pM^{n+1} P^{p-1} H^q- q M^n
 \frac{ \sinh (a P)}{a} P^{p} H^{q-1}, 
\end{equation}
the action of the  generators can be written as
\begin{equation} \begin{split}
   K^m M^n P^p H^q  \prec K=  &
     K^{m+1} M^n P^p H^q - pK^m M^{n+1} P^{p-1} H^q- q K^m M^n\;
            \frac{ \sinh (a P)}{a} P^{p} H^{q-1}, \\[0.2cm]
   K^m M^n P^p H^q  \prec M= &   K^{m} M^{n+1} P^{p} H^{q}, \\[0.2cm]   
   K^m M^n P^p H^q  \prec P= &   K^{m} M^{n} P^{p+1} H^{q}, \\[0.2cm]   
   K^m M^n P^p H^q  \prec H= &   K^{m} M^{n} P^{p} H^{q+1}.   
\end{split} \end{equation}
The endomorphism associated to this action is given by
 \begin{equation} \begin{array}{llll}
       &  \rho(K)=  \displaystyle{\bar{K}- \bar{M} \frac{\partial }{\partial P}-
            \displaystyle \frac{ \sinh (a \bar{P})}{a} \frac{\partial }{\partial H}},&\qquad
        & \rho(M)=   \bar{M}, \\[0.3cm]
         &\rho(P)=  \bar{P},&\qquad
        & \rho(H)=  \bar{H} .
    \end{array}\end{equation}
The adjoint operators respect to the bilinear form (\ref{pairinggext2}),
 \begin{equation} \begin{array}{llll}
         &\rho(K)^\dagger=  \displaystyle{ \frac{\partial }{\partial v}- \bar{x}\; \frac{\partial
}{\partial \mu}- \bar{t}\; \frac{ \sinh (a \frac{\partial }{\partial x})}{a}}, & \qquad
        & \rho(M)^\dagger= \displaystyle{\frac{\partial }{\partial \mu}}, \\[0.3cm]
         &\rho(P)^\dagger=  \displaystyle{\frac{\partial }{\partial x}}, & \qquad
        & \rho(H)^\dagger= \displaystyle{ \frac{\partial }{\partial t}}, 
  \end{array}\end{equation}
allow to obtain the action that characterizes $(U_q(\overline{\mathfrak{g}(1,1)}), \prec,
U_q(\overline{\mathfrak{g}(1,1)}))$
 \begin{equation} \label{lslsls}
  \begin{array}{llll}
     & K\lact f=  \displaystyle{\left[ \frac{\partial }{\partial v} -
                \bar{x}\;\frac{\partial }{\partial \mu}-
                  \bar{t}\;\frac{ \sinh (a \frac{\partial }{\partial x})}{a}\right]f}, &\qquad
     & M\lact f=  \displaystyle{ \frac{\partial }{\partial \mu}f}, \\[0.3cm] 
      &P\lact f=   \displaystyle{\frac{\partial }{\partial x}f}, &\qquad 
     & H\lact f=   \displaystyle{\frac{\partial }{\partial t}f}.  
    \end{array}\end{equation}

Now let study the module $({\cal A},\ract, U_q(\overline{\mathfrak{g}(1,1)}))$, dual of
$(U_q(\overline{\mathfrak{g}(1,1)}), \succ, U_q(\overline{\mathfrak{g}(1,1)}))$ respect to 
the bilinear form (\ref{pairinggext2}). The following lemma will be useful for our task.
\begin{lemma} For every real number  ${s}$ the following relations hold
  \begin{equation} \label{lemama}
     \begin{split}
         M e^{{s} K}= & e^{{s} K} M, \\[0.2cm]
         P e^{{s} K}= &  e^{{s} K}(P- {s} M), \\[0.2cm]
         H e^{{s} K}= &
       e^{{s} K}(H + \frac{\cosh(a(P-{s} M))- \cosh(aP)}{a^2 M}).
     \end{split} \end{equation}
 \end{lemma}
 \begin{proof}[Proof]
Let $\varphi$ be an element  of the  subalgebra generated by  $M,P$ and $H$.
We have
 \begin{equation} \label{conmu}
 \varphi e^{{s} K}= e^{{s} K} e^{{s}\ \rad{K} }(\varphi) ,
 \end{equation}
where
 \begin{equation} \label{adjuntar}
 \rad{a} (a')=a' a - aa'=[a',a], \qquad 
 e^{{-a} } a'  e^{a}=e^{ \rad{a}}(a') ,
 \end{equation}
with $a$ and $a'$ belonging to an associative algebra. Since the subalgebra 
$\langle M, P, H \rangle $ is isomorphic to
$F(\mathbb{R}^3)$ the  system $(M, P, H)$ can be seen as a global chart over $\mathbb{R}^3$ 
(see Ref..~\cite{olmo01}). The derivation $\hat{K}=\rad{K}$ in this chart can be written as
 \begin{equation}
         \hat{K}= -M \frac{\partial}{\partial P}-
   \frac{\sinh(aP)}{a} \frac{\partial}{\partial H}.  
\end{equation}
Denoting by $\phi$ the flow associated to the vector field  $K$   expression
(\ref{conmu}) becomes
 \begin{equation} \label{conmu2}
     \varphi e^{{s} K}=
     e^{{s} K} (\varphi \circ \phi^{s}).
\end{equation}
In order to obtain  $\phi$ it is necessary to compute firstly the integral curves
of $\hat{K}$, which are determined by the autonomous system \cite{olmo02}
 \begin{equation} 
      \dot{m}=  0, \qquad
      \dot{p}= -m, \qquad
      \dot{h}= - \frac{\sinh(ap)}{a} ,
 \end{equation}
whose solution is
\begin{equation}
     m(s)= c_1, \qquad p(s)= c_2 -s c_1,  \qquad
     h(s)= \frac{\cosh(a(c_2- s c_1))}{a^2 c_1}+ c_3.
\end{equation}
From the expression of the integral curves we obtain the flow
\begin{equation}
   \phi^s(m,p,h)= \left( m,\ p-sm,\ h + \frac{\cosh(a(p-sm))- \cosh(ap)}{a^2m}\right).
\end{equation}
Expressions (\ref{lemama}) are obtained
considering successively  $\varphi$
equal to $M,P$ and $H$ in expression (\ref{conmu2}).
\end{proof}

The lemma allows us to  describe the  action of the generators of
$ U_q(\overline{\mathfrak{g}(1,1)})$
over the elements like $e^{{s} K} M^n P^p H^q$. Thus,
\begin{equation} \begin{split}
    K \succ e^{{s} K} M^n P^p H^q = &
        K e^{{s} K} M^n P^p H^q, \\[0.2cm]
    M \succ e^{{s} K} M^n P^p H^q = &
     e^{{s} K} M^{n+1} P^p H^q,  \\[0.2cm]
    P \succ e^{{s} K} M^n P^p H^q = &
                e^{{s} K} (P-{s} M)M^n P^p H^q, \\[0.2cm]
    H \succ e^{{s} K} M^n P^p H^q = &
     e^{{s} K} \left( H + \frac{\cosh(a(P-{s} M))- \cosh(aP)}{a^2 M}\right) M^{n} P^p H^q. 
   \end{split}
\end{equation}
The endomorphism of $ U_q(\overline{\mathfrak{g}(1,1)})$
associated to this action is
\begin{equation} \begin{array}{llll}
      &\lambda(K)=  \bar{K}, &\qquad
     & \lambda(M)=  \bar{M}, \\[0.2cm]
      &\lambda(P)=  \bar{P}- \bar{M} \displaystyle{\frac{\partial }{\partial K}}, &\qquad      
&\lambda(H)= \bar{H}+\displaystyle{\frac{\cosh(a(\bar{P}-\bar{M}  \frac{\partial }{\partial K}))
           - \cosh(a\bar{P})}{a^2  \bar{M}}}. 
   \end{array}
\end{equation}
The adjoint operators respect to the bilinear form
(\ref{pairinggext2}) are given by
\begin{equation}  \begin{array}{llll} 
      &\lambda(K)^\dagger= \displaystyle{\frac{\partial}{\partial v}}, &\qquad
      &\lambda(M)^\dagger= \displaystyle{ \frac{\partial}{\partial \mu}}, \\[0.2cm]
      &\lambda(P)^\dagger= 
          \displaystyle{\frac{\partial}{\partial x}- \bar{v} \frac{\partial }{\partial \mu}}, &\qquad
     & \lambda(H)^\dagger= \displaystyle{\frac{\partial}{\partial t} }+
            \displaystyle{   \frac{\cosh(a(\frac{\partial}{\partial x}-\bar{v} 
               \frac{\partial }{\partial \mu} ))
           - \cosh(a\frac{\partial}{\partial x})}{a^2 \frac{\partial }{\partial \mu}}}. 
   \end{array}
\end{equation}
These expressions lead immediately to the following action  that characterizes the regular module 
$({\cal A},\ract, U_q(\overline{\mathfrak{g}(1,1)}))$
\begin{equation} \label{accdergext}
   \begin{array}{llll}
      & f \ract K=  \displaystyle{\frac{\partial }{\partial v}}f, & \qquad
      & f \ract M=   \displaystyle{\frac{\partial }{\partial \mu}}f, \\[0.3cm]
      & f \ract P=  \displaystyle{ \left[\frac{\partial }{\partial x}-
            \bar{v} \frac{\partial }{\partial \mu} \right]}f, &\qquad
      & f \ract H=  \displaystyle{ \left[ \frac{\partial }{\partial t} +
                \frac{ \cosh (a (\frac{\partial }{\partial x}- \bar{v}
\frac{\partial}{\partial \mu}))-
           \cosh(a \frac{\partial }{\partial x})}{a^2
          \frac{\partial }{\partial \mu}}\right]}f. 
    \end{array}\end{equation}

Now we can state the following theorem that gives a family of local representations induced 
from the abelian subalgebra generated by $M$ and $K$.
\begin{theorem} \label{locales}
Let us consider the character associated to the subalgebra generated by $\{ M,K\}$ defined by
\begin{equation} \label{esta}
             K^m M^n \vdash 1= \alpha^m \beta^n,  \qquad \alpha, \beta \in \mathbb{C}.
\end{equation}
The carrier space, $\mathbb{C}^\uparrow$,
of the  representation of
$U_q(\overline{\mathfrak{g}(1,1)})$
induced by (\ref{esta}) is constituted by the elements of $\cal A$ of the form
 \begin{equation}
e^{\alpha v} e^{\beta \mu} \phi(x,t),
\end{equation}
and is isomorphic to the subspace generated by the elements  $x^p t^q$. The action of the
generators is described by 
\begin{equation} \label{replocgext}
  \begin{split}
 K \vdash \phi(x,t)= & \left[ \alpha- \beta \bar{x}- \bar{t}\; 
\frac{\sinh(a\partial_x)}{a} \right] \phi(x,t) ,\\[0.2cm]        
M \vdash \phi(x,t)= & \beta \phi(x,t) ,\\[0.2cm]
P \vdash \phi(x,t)= & \frac{\partial}{\partial x} \phi(x,t) ,\\[0.2cm]
H \vdash \phi(x,t)= & \frac{\partial}{\partial t} \phi(x,t) .
\end{split}
\end{equation}
\end{theorem}
\begin{proof}[Proof]
The space $\mathbb{C}^\uparrow$ is obtained solving the equivariance  condition
 \begin{equation}
      f \ract K= \alpha f, \qquad f \ract M = \beta f,
 \end{equation}
which according to  (\ref{accdergext}) can be written as
 \begin{equation}\label{equations}
      \frac{\partial}{\partial v} f = \alpha f, \qquad
     \frac{\partial}{\partial \mu} f  = \beta f.
 \end{equation}
The general solution of (\ref{equations}) is, obviously,
 \begin{equation}\label{elementos}
      f= e^{\alpha v} e^{\beta \mu} \phi(x,t),
\end{equation}
with $\phi(x,t)$ a linear combination of monomials
 $x^p t^q$. The  action of  the induced representation (\ref{replocgext}) is obtained by restricting the 
 action (\ref{lslsls}) of the module
$({\cal A}, \lact, U_q(\overline{\mathfrak{g}(1,1)}))$ to the elements (\ref{elementos}).
\end{proof}
\sect{Concluding remarks}\label{conclusions}

Note that in the above computations we have not used at any moment the commutation relations
of the algebra $\cal A$. In other words, we have essentially used the structure of linear
space of  $\cal A$.

The representations characterized by $(\alpha, \beta)$ and $(0, \beta)$ are equivalent. To see 
that it is enough to rescaled $K\to K-\alpha$. So, in the following  we will take $\alpha=0$.

The  action of the Casimir
$C_a= MT +{\displaystyle \frac{1- \cosh(aP)}{a^2}}$
in the induced module
 $(\mathbb{C}^\uparrow, \vdash, U_q(\overline{\mathfrak{g}(1,1)}))$
is given by
\begin{equation}\label{qcasimirequation} 
C_a \vdash \phi(x,t) =\left[ \beta \partial_t+
 \frac{1- \cosh(a \partial_x)}{a^2} \right]
 \phi(x,t).
  \end{equation}
Since the operator $C_a$ is not a constant, the induced  representation is reducible, but the
$q$--Casimir  equation (\ref{qcasimirequation}) gives a condition of irreducibility. For  $\beta \not=0$ 
the generator $M$ acts by means of an invertible  operator and we can  consider 
the central element  $C'_a=M^{-1}C_a$ whose action is
\begin{equation}\label{casimir1} 
C'_a \vdash \phi(x,t) =\left[
 \partial_t+ \frac{1}{\beta}
 \frac{1- \cosh(a \partial_x)}{a^2} \right]
 \phi(x,t).   
\end{equation}

The unitarity of the  representations  requires firstly to fix a  $*$--structure on the
algebra $U_q(\overline{\mathfrak{g}(1,1)})$. We here consider 
a real form that take over the  generators the following value:
\begin{equation}
       K^*= -K, \qquad M^*= -M, \qquad P^*= -P, \qquad H^*= -H.
 \qquad \end{equation}
On the other hand, it is necessary to determine an inner product on the corresponding space.
Note that the action (\ref{replocgext}) is well defined over the functions of 
${\cal H}=L^2(\mathbb{R}^2)$ which are ${\cal C}^\infty$.
For this reason we fix, finally, the commutation relation
$[x,t]=0$
and  make the identification $\mathbb{C}^\uparrow \equiv {\cal H}_\infty$.

Obviously, with the inner product given by
 \begin{equation} 
(\phi \, | \, \phi' )=
         \int_{-\infty}^{+\infty} dx  \int_{-\infty}^{+\infty} dt
            \;  \phi(x,t) \overline{\phi'(x,t)}
 \end{equation}
we have
 \begin{equation} ( P \vdash \phi \, | \, \phi' ) =
      (  \phi \, | \, P^* \vdash \phi' ), \qquad
              ( H \vdash \phi \, | \, \phi' ) =
      (  \phi \, | \, H^* \vdash \phi' ),
    \end{equation}
independently of the values of $\alpha$ and $\beta$
that characterize the  representation from where we induced.
However, to get
 \begin{equation}
     ( M \vdash \phi \, | \, \phi' ) =
      (  \phi \, | \, M^* \vdash \phi' ), \qquad
    ( K \vdash \phi \, | \, \phi' ) =
      (  \phi \, | \, K^* \vdash \phi' )
    \end{equation}
it is necessary to impose the condition
$\bar{\beta}= -\beta$.

In the  limit $a\rightarrow 0$ equation (\ref{casimir1}) reduces to
\begin{equation} \label{casimirgextnodef}
 C'_0 \vdash \phi(x,t) =\left[
 \partial_t- \frac{1}{2\beta}
  \partial_x^2 \right] \phi(x,t),  
\end{equation}
which coincides with the Schr\"odinger equation for a free elemental system  of mass $m$ 
after making $\beta=-im/\hbar$.
Also in the limit $a\rightarrow 0$ the expressions (\ref{replocgext}) become
\begin{equation} \label{replocgextl}
  \begin{array}{llll}
 & K \vdash \phi(x,t)=  \displaystyle{\left( \frac{im}{\hbar} {x}- {t} \partial_x \right) \phi(x,t)} ,
&\qquad   & M \vdash \phi(x,t)=  -\displaystyle{\frac{im}{\hbar}} \phi(x,t) ,\\[0.3cm]
& P \vdash \phi(x,t)=  \displaystyle{\frac{\partial}{\partial x}} \phi(x,t) , &\qquad
& H \vdash \phi(x,t)=  \displaystyle{\frac{\partial}{\partial t} }\phi(x,t) .
\end{array}
\end{equation}

As we mention in the introduction, the motivation for constructing the  $q$--deformed version of the
local representations of $U_q(\overline{\mathfrak{g}(1,1)})$ is  to arrive at a  `deformed'
quantum mechanics, where the deformation parameter of the quantum algebra is introduced as a fundamental
length that can be interpreted as a lattice spacing. This is a first step to construct a $q$--version of
quantum field theory. Results in these directions will be published elsewhere.

\section*{Acknowledgments}
This work has been partially supported by DGES of the Ministerio de Educaci\'on y
Cultura de Espa\~na under Project PB98--0360, and the
Junta de  Castilla y Le\'on (Spain).


\end{document}